\documentclass[11pt]{amsart}

\usepackage[T1]{fontenc}
\usepackage[utf8]{inputenc}
\usepackage{lmodern}
\usepackage{amsmath,amssymb,amsthm}
\usepackage[colorlinks=true,linkcolor=blue!60!black,citecolor=blue!60!black,urlcolor=blue!60!black]{hyperref}
\usepackage{xcolor}
\usepackage{xurl}

\theoremstyle{plain}
\newtheorem{theorem}{Theorem}[section]

\numberwithin{equation}{section}

\DeclareMathOperator{\SL}{SL}
\DeclareMathOperator{\GL}{GL}
\DeclareMathOperator{\Stab}{Stab}
\DeclareMathOperator{\Sym}{Sym}
\DeclareMathOperator{\Char}{Char}

\newcommand{\GG}{\mathbb{G}}

\renewcommand{\AA}{\mathbb{A}}

\begin{document}

\title[Weitzenb\"ock's theorem in arbitrary characteristic]{Weitzenb\"ock's theorem in arbitrary characteristic}

\author{Dominic Bunnett}
\address{Institut f\"ur Mathematik, Technische Universit\"at Berlin, Stra\ss e des 17.~Juni 136, 10623 Berlin, Germany}
\email{bunnett@tu-berlin.de}


\date{}

\begin{abstract}
Let \(k\) be a field and \(V\) a finite-dimensional rational representation of the additive group \(\GG_{a,k}\).
We show that the ring of invariants \(k[V]^{\GG_a}\) is finitely generated in every characteristic, extending Weitzenb\"ock's classical theorem beyond characteristic zero.
\end{abstract}

\maketitle

\section{Introduction}

Let \(k\) be a field and consider a representation \(\rho : \GG_{a,k} \to \GL(V)\) and write \(k[V]=\operatorname{Sym}_k(V^\vee)\).
When \(\Char k =0\), Weitzenb\"ock's Theorem says that the ring of invariants \(k[V]^{\GG_a}\) is finitely generated \cite{Weitzenbock}.
This note proves that the same finite-generation statement holds in positive characteristic.

\begin{theorem}\label{thm:main}
Let \(k\) be a field and let \(V\) be a finite-dimensional rational
representation of \(\mathbb G_{a,k}\). Then \(k[V]^{\mathbb G_{a,k}}\) is a
finitely generated \(k\)-algebra.
\end{theorem}

The question arose for the author in studying non-reductive geometric invariant theory in positive characteristic.
Modin extended the \(\widehat U\)-theorem over affine Noetherian bases of arbitrary characteristic, obtaining geometric quotients on specified loci for graded unipotent actions under suitable stabiliser hypotheses \cite{Modin}.
Cooper and Modin subsequently developed an intrinsic version that constructs good moduli spaces after rigidification under hypotheses on the centre and relative inertia \cite{CooperModin}.
The theorem above answers the affine invariant-theoretic question for a single additive group: the affine categorical quotient of every finite-dimensional linear action is of finite type.

Seshadri proved the characteristic-zero theorem by reduction to \(\mathrm{SL}_2\) \cite{Seshadri}:
Let \(G = \SL_{2,k}\) and \(U \cong \GG_{a,k}\) be the unipotent upper triangular matrices.
We take \(k\) algebraically closed (we discuss this below in the actual proof).
Consider the homogeneous fibre space \(G \times^{U}V\) which is the quotient of \(G \times V\) by the \(U\)-action defined by \(u \cdot (g,v) = (gu^{-1} , uv)\).
As in Grosshans' observable subgroup theory \cite[Section 9]{Grosshans}, we can identify 
\[k[G\times^{U}V]^{G} \cong k[V]^{U},\]
where the ring on the left is the invariant global sections of the quasi-affine variety \(G \times^U V\).
See also \cite{DoranKirwan}.

The representation \(\rho\) extends to a \(G\)-representation \(\hat{\rho} : G \to \GL(V)\) in characteristic zero.
Thus we get a \(G\)-equivariant isomorphism
\[G \times^{U}V \cong \left(G/U \right)\times V,\quad[g,v]\mapsto (gU,gv).\]
Then as \(G/ U \cong \AA^2\setminus\{0\}\), we can extend uniquely the \(G\)-action to \(\AA^2\).
Thus
\[k[\AA^2 \times V]^{G} = k[G\times^{U}V]^{G} \cong k[V]^{U}\]
since \(\AA^2\) is normal and \(\{0\}\) has codimension 2.
We know the first ring is finitely generated since \(\SL_2\) is reductive and we are done.

In positive characteristic, however, not every linear \(\GG_a\)-action extends to \(\SL_2\).
Indeed the terrain looks very different.
Fauntleroy proved finite generation over an algebraically closed field when the fixed subspace \(V^{\mathbb G_a}\) is a hyperplane \cite{Fauntleroy}.
He uses common multiples of additive polynomials to construct a polynomial subalgebra over which the invariant ring is finite.
In the paper he gives an example showing that without this assumption, one cannot use the same method.

Over an algebraically closed field of characteristic two, Tanimoto proved that the invariant ring of every nontrivial four-dimensional representation is a polynomial ring in three variables \cite{Tanimoto}.

Roberts, Freudenburg, and Daigle–Freudenburg constructed nonlinear \(\mathbb G_a\)-actions with non-finitely generated invariant rings in characteristic zero \cite{Roberts,Freudenburg,DaigleFreudenburg}.
Dufresne and Maurischat \cite{DufresneMaurischat} constructed characteristic-free versions of these three examples and proved that their invariant rings are finitely generated in every positive characteristic.

Examples of nonlinear \(\GG_a\)-actions on affine spaces with non-finitely-generated invariant rings in every positive characteristic follow by combining Doran–Giansiracusa \cite[Corollary 1.3]{DoranGiansiracusa} with Castravet–Laface–Tevelev–Ugaglia  \cite[Theorem 1.2]{CastravetLafaceTevelevUgaglia}.
Doran and Giansiracusa identify \(\operatorname{Cox}(\overline M_{0,n})\) with the invariant ring of a nonlinear \(\mathbb G_a\)-action on affine space.
For \(n\ge 10\), this Cox ring is not finitely generated in any positive characteristic, by the non-polyhedrality theorem of Castravet, Laface, Tevelev, and Ugaglia \cite[Theorem 1.2]{CastravetLafaceTevelevUgaglia}.

Maguire has proposed a counterexample to finite generation in the linear case \cite{Maguire}.
For the action intended in \cite[Thm.~5.1]{Maguire}, the proposed identification in Theorem 3.3 already fails in degree two of the filtration grading: divisibility by \(x_1^2\) does not imply membership in the second symbolic power.
It was in reviewing the proposed counterexample \cite{Maguire}, which is a very nice and geometric idea, that the author found the proof presented here.

In this note, we show that in fact the same proof as above does work with a minor tweak.
The slight difference in perspective is that the homogeneous fibre space \(G \times^U V\) is still a vector bundle over \(\AA^2 \setminus \{0\}\), but the trivialisation of the extension of the bundle to \(\AA^2\) may not be \(\SL_2\)-equivariant, but this does not turn out to matter.

\section{Proof of the main theorem}

\begin{proof}[Proof of Theorem~\ref{thm:main}]
First we note that we may assume that \(k\) is algebraically closed, as taking invariants commutes with flat base change and finite generation descends along \(k \subset \bar{k}\).

Set \(G=\SL_{2,k}\), \(X=\mathbb A^2_k\) and \(X^\circ=X\setminus\{0\}\).
The first-column map \(G\to X^\circ\) is a torsor under \(U=\Stab_G(e_1)\cong\mathbb G_{a,k}\).
Indeed, on \(D(x)\) and \(D(y)\) one can complete the first column \((x,y)^{\mathsf T}\) by \((0,1/x)^{\mathsf T}\) and \((-1/y,0)^{\mathsf T}\), respectively.
Write \(s_x,s_y\) for these sections.
Over either chart, write \(g=s_i(z)u\).
In the coordinates \((z,u,w)\) given by \(w=uv\), the action \(h\cdot(g,v)=(gh^{-1},hv)\) takes \((z,u,w)\) to \((z,uh^{-1},w)\).
Thus the quotient on each chart is the product of that chart with \(V\).
On \(D(xy)\) the fibre coordinates satisfy \(w_x=(s_x^{-1}s_y)w_y\), so these local quotients glue by a regular linear change of coordinates to a vector bundle \(G\times^U V\to X^\circ\).
Left multiplication on \(G\) induces its \(G\)-equivariant structure.
Let \(\mathcal W^\circ\) denote the locally free sheaf of sections of the dual bundle, equipped with the induced \(G\)-linearisation.

Using Hartogs' theorem for reflexive sheaves we extend \(\mathcal W^\circ\) and the \(G\)-linearisation to a \(G\)-equivariant reflexive sheaf on \(X\) which turns out to be locally free.
Explicitly, for \(j\colon X^\circ\hookrightarrow X\), the extension \(\mathcal W=j_*\mathcal W^\circ\) is coherent reflexive by the depth-two extension theorem and, moreover, it is locally free because \(X\) is a regular surface \cite[Tags \href{https://stacks.math.columbia.edu/tag/0EBJ}{0EBJ} and \href{https://stacks.math.columbia.edu/tag/0B3N}{0B3N}]{Stacks}.

Let \(a,p\colon G\times X\to X\) be the action and projection.
On \(G\times X^\circ\), the given linearisation is an isomorphism \(a^*\mathcal W^\circ\xrightarrow{\sim}p^*\mathcal W^\circ\).
This isomorphism extends uniquely as sections of the corresponding locally free Hom sheaves across the codimension-two subset \(G\times\{0\}\).
The extension is an isomorphism \(a^*\mathcal W\xrightarrow{\sim}p^*\mathcal W\) by uniqueness of reflexive extensions.
The identity and cocycle laws follow by the same uniqueness argument on \(X\) and \(G\times G\times X\), respectively.

Let
\[
Y=\operatorname{Spec}_X(\Sym_{\mathcal O_X}\mathcal W),
\qquad
Y^\circ=Y\times_X X^\circ.
\]
The morphism \(\pi\colon Y\to X\) is a vector bundle, affine and of finite type.
Since \(X\) is affine, \(Y\) is affine of finite type over \(k\).
Over \(X^\circ\), the bundle \(\mathcal W^\circ\) is associated to \(V^\vee\), so its relative spectrum \(Y^\circ\) is \(G\times^U V\).
In particular, if \(\pi^\circ\colon Y^\circ\to X^\circ\) is the projection, the degree-\(d\) part of \(\pi^\circ_*\mathcal O_{Y^\circ}=\Sym_{\mathcal O_{X^\circ}}\mathcal W^\circ\) is \(\Sym^d\mathcal W^\circ\).
This is the grading by polynomial degree along the fibres.

The complement \(Y\setminus Y^\circ\) is the fibre over \(0\) and has codimension two in the smooth variety \(Y\).
By Hartogs, restriction therefore gives an isomorphism \(k[Y]\xrightarrow{\sim}\Gamma(Y^\circ,\mathcal O_{Y^\circ})\).
The restriction map is \(G\)-equivariant, so it induces an isomorphism on \(G\)-invariants.

Let us recall explicitly the transfer principle of Grosshans \cite[Section 9]{Grosshans}.
Restriction to the fibre over \(e_1\), via \(v\mapsto[1,v]\), takes
\(G\)-invariant regular functions on \(Y^\circ\) to \(U\)-invariant
polynomials on \(V\), since \(U\) is the stabiliser of \(e_1\).
Conversely, for \(h\in k[V]^U\), the formula \([g,v]\mapsto h(v)\)
defines a \(G\)-invariant regular function on \(Y^\circ\).
Indeed, in the bundle trivialisations over \(D(x)\) and \(D(y)\)
it is the polynomial \(h\) on the fibre. These local functions agree
on the overlap since the transition functions act through \(U\).
The two constructions are inverse since \([g,v]=g\cdot[1,v]\).
Thus the two restriction maps give isomorphisms of \(k\)-algebras
\[
k[Y]^G
\xrightarrow{\sim}\Gamma(Y^\circ,\mathcal O_{Y^\circ})^G
\xrightarrow{\sim}k[V]^U.
\]
Geometric reductivity and Nagata's finite-generation theorem
\cite{Haboush,Nagata} give finite generation of \(k[Y]^G\), and hence
of \(k[V]^U\).
\end{proof}

\section*{Acknowledgements and use of AI}

I would like to thank Vicky Hoskins for feedback on a first draft which greatly helped with clarity of the exposition.

As part of good practice, I detail what part AI played in this.
I used ChatGPT Sol 6.2 and Qwen3.5:9B as `research assistants'.
Qwen is my own personal model (an open source LLM which runs on my laptop) which has been tweaked and which I have done some fine tuning training to.
It is capable as a coding agent and not capable of doing research level math (yet).
It can manage graduate level algebraic geometry, not more.
My impression is however, that it is not entirely an issue of the model, but of the RAM I have on my machine.
There are details on my website: \emph{\href{https://dombunnett.github.io/trough/}{the trough}}.

Specifically, I found an error in the nefness and extremality argument in the proof of \cite[Thm.~5.1]{Maguire}.
ChatGPT confirmed this and then found the filtration error described above.
Then I used Qwen to help me run lots of software experiments looking for finite generating sets in the example, which it found.
Then using both models, we started to come up with proofs for entire families of \(\GG_a\)-actions which worked very well, but there was not a visible path to the general case.
At some point when looking at one particular family, ChatGPT told me that (in the notation of the proof above) \(\mathcal W^\circ\) was a `vector bundle' on \(\AA^2 \setminus \{0\}\) (a completely trivial fact, but one I had never thought of before).
At this point, the proof was then immediately clear to me.

As far as writing is concerned, all the mathematics was written by me and of course I take full responsibility for content.



\end{document}